\documentclass[12pt]{amsart} 

\usepackage[psamsfonts]{amssymb} 
\usepackage{amsfonts,amsmath}

\newtheorem{thmA}{Theorem}

\newtheorem{corA}[thmA]{Corollary}

\newtheorem{theorem}{Theorem}[section]
\newtheorem{prop}[theorem]{Proposition}
\newtheorem{proposition}[theorem]{Proposition}
\newtheorem{lemma}[theorem]{Lemma}
\newtheorem{corollary}[theorem] {Corollary}

\theoremstyle{remark}
\newtheorem{remark}[theorem]{Remark}

\theoremstyle{definition}

\def\Z{\mathbb Z}

\def\G{\Gamma}
\def\g{\gamma}

\def\autn{{\rm{Aut}}(F_n)}
\def\aut{{\rm{Aut}}}
\def\out{{\rm{Out}}}
\def\outn{{\rm{Out}}(F_n)}
\def\gln{{\rm{GL}}(n, \Z)}
\def\glm{{\rm{GL}}(m, \Z)}
\def\ian{{\rm{IA}_n}}
\def\ion{{\overline{\rm{IA}}_n}}
\def\L{\Lambda}

\def\hom{{\rm{Hom}}}
\def\qh{{\rm{QH}}}
\def\ad{{\rm{ad}}}
\def\<{\langle}
\def\>{\rangle}

\title{Actions of higher-rank lattices on free groups}
\author{Martin R. Bridson and Richard D. Wade}   
\begin{document}
\begin{abstract} If $G$ is a semisimple Lie group of real rank at least 2
and $\Gamma$ is an irreducible lattice in $G$, then every homomorphism from $\G$ to the outer automorphism
group of a finitely generated free group has finite image.
\end{abstract}


\address{Mathematical Institute, 24-29 St Giles', Oxford OX1 3LB, UK}
\email{bridson@maths.ox.ac.uk}
\email{wade@maths.ox.ac.uk}

\thanks{This work is funded by the EPSRC of Great Britain: Bridson is supported by a Senior Fellowship,
Wade is supported by a DTA grant.} 

\maketitle
\section{Introduction}

In recent years, a powerful body of mathematics has emerged from
efforts to extend rigidity phenomena from the context of irreducible lattices in
semisimple Lie groups to a wider context that embraces mapping class groups of surfaces
and automorphism groups of finitely generated
free groups (see \cite{Z} and \cite{BV} for references). An important focus
of these efforts has been the conjecture that every map from 
an irreducible, higher-rank lattice $\Gamma$
to a mapping class group or the automorphism group of a finitely generated free group must have finite image. This was proved in the case of mapping class groups by Farb, Kaimanovich and Masur \cite{FM, KM}; subsequent
 proofs of their result have elucidated different aspects of the
geometry of mapping class groups and their subgroups \cite{BFuji, DMS, DGO}. 

When $\Gamma$ is non-uniform, one obtains a short proof in the mapping class group case 
by combining the Normal Subgroup Theorem \cite[Theorem 8.1.2]{nsgp}
with the fact that all solvable subgroups of 
mapping class groups are virtually abelian \cite{ivanov, birman-lub-McC}. 
A similar argument,
using \cite{alibeg} and \cite{bfh} in place of  
\cite{ivanov} and \cite{birman-lub-McC}, shows that any homomorphism from
$\Gamma$ to the outer automorphism group of a finitely generated free group must also factor through a finite group; see \cite{bridson-farb}. Our main objective in this article is the 
corresponding result for uniform lattices.

\begin{thmA}\label{thm} Let $\G$ be a group. Suppose that no subgroup
of finite index in $\G$ has a
normal subgroup that maps surjectively to $\Z$.
Then every homomorphism from $\G$ 
to the outer automorphism group of a finitely generated free group
has finite image.
\end{thmA}

In Theorem \ref{thm} and the variations on it in Section~\ref{s:refinements}, we do not assume that $\G$
is finitely generated.  

We say that a group satisfying the hypothesis of Theorem~\ref{thm} is \emph{$\mathbb{Z}$--averse}.
The Normal Subgroup Theorem of Kazhdan and Margulis \cite{nsgp} tells us that
irreducible lattices in connected, higher-rank, semisimple Lie groups with finite centre have no infinite normal subgroups
of infinite index. Since such lattices are not virtually cyclic, it follows that they are $\mathbb{Z}$--averse.

\begin{corA}\label{c:main} 
If $G$ is a connected, semisimple Lie group of real rank at least 2 that has finite centre,
and $\Gamma$ is an irreducible lattice in $G$, then every homomorphism from $\G$ to the outer automorphism
group of a finitely generated free group has finite image.
\end{corA}

An additional argument allows one to remove the hypothesis that $G$ has finite centre (see Remark \ref{r:ext}). Further examples of $\mathbb{Z}$--averse groups come from  Bader and Shalom's recent work 
on the Normal Subgroup Theorem \cite{BS}.
If a hereditarily just infinite group is not virtually cyclic then it  is  $\mathbb{Z}$--averse; examples are described in \cite{Wilson}.

Corollary \ref{c:main} has implications for the Zimmer programme \cite{Z}. Specifically, it allows one to extend Farb and Shalen's theorem about actions of higher-rank lattices on 3--manifolds to the general case, removing the non-uniform hypothesis from Theorem II of \cite{FS} and improving their Theorem III as follows:

\begin{theorem}
Let $\G$ be an irreducible uniform lattice in a semisimple Lie group of real rank at least 2, and let $M$ be any closed, oreintable, connected 3--manifold. Then for every action $\G \to {\rm{Homeo}}(M)$, the image of $\G$ in $\aut(H_*(M,\mathbb{Q}))$ is finite. \end{theorem}

Our proof of Theorem \ref{thm} relies heavily on recent results of Bestvina and Feighn \cite{BF}, Dahmani, Guirardel
and Osin \cite{DGO}, and Handel and Mosher \cite{HM}. The work of 
Bestvina and Feighn was inspired in part by the desire to prove Corollary
\ref{c:main}, following the lines of the proof given in the case of mapping
class groups by Bestvina and Fujiwara \cite{BFuji}, which invokes
Burger and Monod's theorem that irreducible
 lattices in higher-rank Lie groups
have trivial bounded cohomology \cite{BM}. One can replace this
use of bounded cohomology with an argument of
Dahmani, Guirardel and Osin that applies small cancellation 
theory to the study of purely pseudo-Anosov subgroups;
this is used in  \cite{DGO} to prove an analogue of Theorem \ref{thm} for homomorphisms
to mapping class groups.

Whether one uses bounded cohomology or the alternative endgame from
\cite{DGO}, the key step in the Bestvina--Feighn--Fujiwara approach
is to get a finitely generated
subgroup of $\outn$ to act in a suitable way on a hyperbolic metric space.
Bestvina and Feighn \cite{BF} construct such actions for
subgroups of $\outn$ that 
contain a fully irreducible automorphism, and hence deduce that
a higher-rank lattice cannot map onto such a subgroup.   (Hamenst\"adt \cite{ursula} has a different way
of constructing bounded cohomology classes for these subgroups.)  
 If one could 
construct suitable actions for more general 
subgroups of $\outn$, then Theorem \ref{thm} would follow. A significant step in this direction was taken recently by Handel and Mosher \cite{HM}, who proved that if
a subgroup $H<\outn$ does not contain the class of a
fully irreducible automorphism, then $H$ has a subgroup of finite index that
leaves the conjugacy class of a proper free factor of $F$ invariant. 
Handel and Mosher also indicate that they hope to extend their
work so as to prove Corollary \ref{c:main}
along the lines sketched above.

Arguments of a quite different sort allow one to see that a homomorphism
from a uniform higher-rank lattice $\G$ to $\autn$ cannot contain any polynomially
growing automorphisms of infinite order: on the one
hand, Piggott \cite{adam} proves that the homomorphism $\autn\to\glm$ given by the action
of $\autn$ on the first homology of some characteristic subgroup of finite
index in $F_n$ will map
a power of such an automorphism  to a 
non-trivial unipotent; on the other hand, Margulis superrigidity
implies that the image of any homomorphism $\G\to\glm$ can contain only semisimple elements.

Our proof of Theorem \ref{thm} proceeds as follows.
In Proposition
\ref{p:IA_n} we shall use the results of Handel--Mosher and Dahmani--Guirardel--Osin to see that if $\G$ is $\mathbb{Z}$--averse, then
the image of every homomorphism 
$\Gamma\to\outn$ will have a subgroup of finite index that lies in the
kernel $\ion$ of the map $\outn\to\gln$ given by the action of $\outn$ on the
first homology of $F_n$. In Section~\ref{s:BL} we use Lie methods, \`a la Magnus, to prove
that every  non-trivial subgroup of $\ion$ maps onto $\Z$ (Corollary \ref{c:BL}). (This
result also appears in Bass and Lubotzky's work on central series \cite{BL}.)
As no finite index subgroup of $\G$ maps onto $\Z$, this completes the
proof of Theorem~\ref{thm}.

We thank Laurent Bartholdi for a helpful conversation concerning central filtrations.

\section{Proofs} \label{s:proofs}

We fix a $\mathbb{Z}$--averse group $\G$. We do not assume that $\G$ is finitely generated.

\subsection{Controlling the action of $\G$ on homology}

\begin{proposition}\label{p:IA_n} For every subgroup
of finite index $\Lambda\subset\Gamma$ and every homomorphism
$\phi:\Lambda\to\outn$, the intersection $\phi(\Lambda)\cap\ion$
has finite index in $\phi(\Lambda)$.
\end{proposition} 

\begin{proof} The proof is by induction on $n$. The case $n=1$ is trivial.
${\rm{Out}}(F_2)$ has a free subgroup of finite index and  no subgroup
of finite index in $\G$ can map onto a free group, so every 
homomorphism $\Lambda\to{\rm{Out}}(F_2)$ has finite image.

Suppose $n\ge 3$. Recall that $\psi\in\autn$ (and its image in
$\outn$) is said to be \emph{fully
irreducible} if no power of $\psi$ sends a 
proper free factor of $F_n$ to a conjugate of itself.  Let $[\psi]$ denote the image of $\psi$ in $\outn$.
 Using the actions constructed in \cite{BF} and drawing on the
approach to small cancellation theory developed in \cite{DelzG}, Dahmani, Guirardel and Osin \cite{DGO} prove that if
$\psi$ is fully irreducible then for some positive integer $N$, the normal
closure of $[\psi]^N$ is a free group. It follows that any subgroup
of $\outn$ that contains a fully irreducible automorphism also contains an
infinite normal subgroup that is free. In particular, $\phi(\Lambda)$ cannot
contain a fully irreducible automorphism.

According to \cite[Theorem 1.1]{HM}, if $\phi(\Lambda)$ does not contain a fully irreducible
automorphism then  a subgroup of finite
index $H \subset \phi(\Lambda)$ leaves a free factor of $F_n$ invariant up to conjugacy;
say $F_n=L\ast L'$, where $\psi(L) = g_\psi^{-1}Lg_\psi$ for all
$[\psi]\in H$. Note that the image in $\out(L)$ of $x\mapsto g_\psi \psi(x) g_\psi^{-1}$,
which we denote $[\psi]_{L}$, depends only on the image
of $\psi$ in $\outn$, and that $[\psi]\mapsto [\psi]_L$ defines a 
homomorphism from $H$ to $\out(L)$. Likewise, the  action
on the quotient  $F_n/\<\!\<L\>\!\>$  induces a homomorphism $H\to {\out}(L')$.
By induction, we know that the induced action of $H$ on the abelianization
 of both $L$ and $L'$ factors through a finite group. Thus the
action of $H$ on the abelianization of $F_n=L\ast L'$ lies in a
block triangular subgroup (with
respect to a basis that is the union of bases for $L$ and $L'$)
$$
\begin{pmatrix}
 G & 0\\
\ast & G'
\end{pmatrix} \leq \gln
$$
where $G$
and $G'$ are finite. This matrix 
group is finitely generated and virtually abelian,
whereas  $\Gamma,$ and therefore $H$, does not have a subgroup of finite index that maps
onto $\Z$. Thus the action of $H$ on the homology of $F_n$ factors through a finite group, and hence that of $\phi(\Lambda)$ does too, i.e.~$\phi(\Lambda)\cap\ion$
has finite index in $\phi(\Lambda)$. This completes the induction.
\end{proof}

\subsection{Central filtrations of $\ian$ and $\ion$.} \label{s:BL}

Let $\gamma_c$ be the $c$th term in the lower central series of $F_n$; so $\gamma_1=F_n$ and $\gamma_c=[\gamma_{c-1},F_n]$.  Let $\Gamma_c=F_n/\gamma_c.$ As $\gamma_c$ is characteristic, there is a natural map $\autn \to \aut(\Gamma_c).$  Let $G_{c-1}$ be the kernel of this map.  Note that $G_0=\autn$ and $G_1=\ian$. Magnus \cite{magnus} showed that $\cap_{i=1}^\infty \gamma_c=0;$ it follows that $\cap_{i=0}^\infty G_c=0.$  In fact, 
Andreadakis proved that $G_1, G_2, G_3, \ldots$ forms a central series for $G_1$ with each quotient  a finitely generated free abelian group \cite{andreadakis}. It now seems natural to regard his result in the context of \emph{higher Johnson homomorphisms} \cite{day,M}. We include a proof for the convenience of the reader. Let $L_{c}=\gamma_c/\gamma_{c+1}.$ Our commutator convention
is $[x,y]:=x^{-1}y^{-1}xy$.

\begin{prop} \label{p:andreadakis}
 For  $c \geq 1$ there exists a homomorphism \[ \tau_c:G_c \to \hom(H_1(F_n),L_{c+1}) \] such that $\ker(\tau_c)=G_{c+1}.$
\end{prop}

\begin{proof}
For all $\psi \in G_c$ and $x \in F_n$ there exists $w_x \in \g_{c+1}$ with $\psi(x)=xw_x.$ Define $\tau_c(\psi)([x])=w_x\gamma_{c+2}$.  Note that if $x \in \g_2=[F_n,F_n]$ then\begin{equation} \label{e:mod}
\psi(x)\g_{c+2}=x\g_{c+2}. \end{equation}
Indeed the commutator relations
$[x,yz]=[x,z][z,[y,x]][x,y]$ and $[xy,z]=[x,z][[x,z],y][y,z]$ imply \begin{align*} \psi([x,y])&=[\psi(x),\psi(y)] \\
&=[xw_x,yw_y] \\
&=[xw_x,w_y][w_y,[y,xw_x]][xw_x,y] \\
&=[xw_x,w_y][w_y,[y,xw_x]][x,y][[x,y],w_x][w_x,y]
\end{align*} and $[xw_x,w_y],[w_y,[y,xw_x]],[[x,y],w_x],[w_x,y]$ all lie in $\gamma_{c+2}$. It follows easily that $\tau_c$ is well-defined and a homomorphism. The automorphism $\psi$ belongs to $G_{c+1}$ if and only if $w_x$ is contained in $\gamma_{c+2}$ for all $x$ in $F_n$, hence $\ker\tau_c=G_{c+1}$. \end{proof}

Proposition \ref{p:andreadakis} shows that $G_c/G_{c+1}$ is isomorphic to a subgroup of the group $\hom(H_1(F_n),L_{c+1})$. The abelian group $L_{c+1}$ is a subgroup of the finitely generated free nilpotent group $\Gamma_{c+2}$; thus it is finitely generated and free abelian.

\begin{corollary} \label{c:ian}
If $c \geq 1$ then $G_c/G_{c+1}$ is a finitely generated free abelian group.
\end{corollary}

If $S$ is a non-trivial subgroup of $\ian=G_1$ then there exists $c \geq 1$ such that $S  \leq G_c$ and $S\not \leq G_{c+1}.$

\begin{corollary} \label{c:abelianisation ian}
Every non-trivial subgroup of $\ian$ maps onto $\Z$.
\end{corollary}

We would like to extend this analysis to subgroups of $\ion$.  Let $H_c$ be the image of $G_c$ under the projection $\pi:\autn\to \outn$.  Our goal in the remainder of this section is to prove the analogue of Corollary \ref{c:ian} for $H_1,H_2,\ldots$
(cf.~\cite{BL}).  We make use of a theorem of Magnus that $L=\oplus_{c=1}^\infty L_c,$ along with the bracket operation induced by commutation in $F_n,$ is a \emph{free Lie $\Z$--algebra} generated by the images in $H_1(F_n)=L_1$
of a basis for the free group.  This theorem and the required background on free Lie algebras is explained in Chapter 5
of  \cite{MKS}.  Let $p$ be a positive integer and let  $(L)_p$ be the free Lie algebra obtained by taking the tensor product of $\Z/p\Z$ with $L$. Let $L^c$ be the quotient algebra $L/\oplus_{i>c}L_i.$ We will need the following fact, whose
proof is sketched in Exercise 3.3 of Chapter 2 in \cite{Bour}. Throughout, $Z(A)$ denotes the centre of a Lie algebra $A$.

\begin{prop} \label{p:centres}
 If $n \geq 2,$ then $Z(L)=Z((L)_p)=0$ and $Z(L^c)$ is the image of $L_c$ under the quotient map
$L\to L^c$.
\end{prop}

As ${\overline{\rm{IA}}_1}$ is trivial, we restrict ourselves to $n \geq 2$ for the remainder of this section. 
 For each non-trivial $y \in F_n$ there exists a unique $c$ such that $y \in \g_c$ and $y \not \in \g_{c+1}.$ We  identify $y$ with its (non-trivial) image in the submodule $L_c$ of the Lie algebra $L$.
Let $\ad:F_n \to \autn$ be the map induced by the action of $F_n$ on itself by conjugation.

\begin{corollary} \label{c:centres}
The kernel of the map $F_n \to \aut(\Gamma_{c+1})$ induced by $\ad$ is $\g_c.$  Hence $\ad(y)$ belongs to $G_c$ if and only if $y$ is in $\g_{c}$
\end{corollary}

\begin{proof} If $y\in\g_c$ then $x^{-1}y^{-1}xy$ lies in $\gamma_{c+1}$ for each $x \in F_n$, therefore
 $\ad(y)$ lies in $G_c$. Conversely, if $y$ is not in $\g_c$ then by Proposition \ref{p:centres},
its image in $L^c$ is not central. As $L^c$ is generated by the images of a basis $x_1,\dots,x_n$ for $F_n$, this means 
that some $[x_i,y]\neq 0$ in $L^c$, hence $x_i^{-1}y^{-1}x_iy_i\notin \g_{c+1}$. 
\end{proof}

\begin{theorem}\label{t:allfree}
For $c \geq 1$ the sequence $$0 \to \g_c/\g_{c+1} \to G_c/G_{c+1} \to H_c/H_{c+1} \to 0$$ is an exact sequence of free abelian groups, where the second and third maps are induced by $\ad$ and $\pi$ respectively.
\end{theorem}

\begin{proof}
Surjectivity of the map $G_c/G_{c+1} \to H_c/H_{c+1}$ is trivial, and exactness of the remaining maps follows from Corollary~\ref{c:centres}.  It remains to show that $H_c/H_{c+1} \cong (G_c/G_{c+1})/\ad(\g_c)$ is  torsion free. Suppose that $\psi \in G_c$ and there exists $y \in \g_c$ and $p\ge 1$
such that $\ad(y)G_{c+1}=\psi^p G_{c+1}.$  For each $x_i$ in a free generating set of $F_n$ we have $$\psi(x_i)=x_iw_i$$ for some $w_i$ in $\g_{c+1}.$  Equation (\ref{e:mod}) in Proposition~\ref{p:andreadakis} tells us 
that $\psi(w_i)\g_{c+2}=w_i\g_{c+2}$, therefore $$\psi(x_i)^p\g_{c+2}=x_iw_i^p\g_{c+2}.$$
As $\ad(y)G_{c+1}=\psi^p G_{c+1},$ this rearranges to $$x_i^{-1}y^{-1}x_iy\g_{c+2}=w_i^p\g_{c+2},$$ therefore $[x_i,y]=0$ in the associated Lie algebra $(L)_p.$ However, $(L)_p$ is generated by the images of $x_1,\ldots,x_n,$ so the image of $y$ lies in $Z((L)_p)=0.$ It follows that $y$ lies in the kernel $pL$ of the map $L \to (L)_p.$  Hence there exists $y_0$ in $\g_{c}$ such that $y\g_{c+1}=y_0^p\g_{c+1},$  and $\ad(y_0)^pG_{c+1}=\psi^pG_{c+1}.$  But $G_c/G_{c+1}$ is free abelian, so $\ad(y_0)G_{c+1}=\psi G_{c+1}.$
\end{proof}

The final fact we need concerning $\outn$ is well known.

\begin{lemma}
The intersection $\cap_{i=1}^\infty H_i$ is trivial in $\outn$.
\end{lemma}

\begin{proof} 
By Proposition I.4.9 of \cite{LS}, if $\psi(u)$ is conjugate to $u$ in $\G_c$ for all $c,$ then $\psi(u)$ is conjugate to $u$ in $F_n.$  And if $\psi$ takes every element of $F_n$ to a conjugate of itself, then $\psi$ is inner
(\cite{G}, Lemma 1). 
\end{proof}

It follows
that if  $S$ is a non-trivial subgroup of $\ion=H_1$ then there exists $c \geq 1$ such that $S  \leq H_c$ and $S \not \leq H_{c+1}.$
In Theorem \ref{t:allfree}
we saw that for $c \geq 1$ the quotient $H_c/H_{c+1}$ is a finitely generated free abelian group.

\begin{corollary} \label{c:BL}
Every non-trivial subgroup of $\ion$ maps onto $\Z$.
\end{corollary}

As we explained at the end of the introduction, this completes the proof of Theorem \ref{thm}.

\section{Alternative Hypotheses} \label{s:refinements}

It emerges from the proofs in the previous section that one
can weaken the hypotheses of Theorem \ref{thm} as follows. Note that
since we have not assumed $\G$ to be finitely generated, condition (1)
is not equivalent to assuming that every finite index subgroup of $\Gamma$ has finite abelianization.

\begin{theorem} \label{t:refinement} Let $\G$ be a group. Suppose that each finite-index
subgroup $\Lambda\subset\G$ satisfies the following conditions:
\begin{enumerate} \item $\Lambda$ does not map surjectively to $\Z$;
\item $\Lambda$ does not have a quotient containing a non-abelian, normal, free subgroup.
\end{enumerate}
Then every homomorphism $\G\to\outn$ has finite image.
\end{theorem}

\begin{proof}
The only additional argument that is needed concerns the normal closure $I$ of $[\psi]^N$ in $\outn$,
as considered in the second paragraph of the proof of Proposition~\ref{p:IA_n}. We must exclude
the possibility that the intersection of $I$ with the image of $\phi:\Lambda \to \outn$ 
is cyclic, generated by $[\psi]^m$ say.  But if this were the case, 
 $\langle [\psi]^m\rangle$ would be normal
in $\phi(\Lambda)$.
Since the normalizer in $\outn$ of the subgroup generated
by any fully irreducible element is virtually cyclic \cite{bfh1},  it would follow that $\phi(\Lambda)$ itself was virtually cyclic, contradicting the fact that no finite-index subgroup of $\Lambda$ maps onto $\Z.$
\end{proof}

\begin{remark}
The class of groups that satisfy the hypotheses of Theorem \ref{t:refinement} is closed under 
certain extension operations. For example, if  $1 \to A \to \hat{\G} \to \G \to 1$ is a short exact sequence and,
in the notation of Theorem \ref{t:refinement}, we suppose that
\begin{itemize}
 \item every finite-index subgroup of $A$ satisfies (2),
\item every finite-index subgroup of $\hat{\G}$ satisfies (1), and
\item every finite-index subgroup of $\G$ satisfies both (1) and (2),
\end{itemize}
then an elementary argument shows that every finite-index subgroup of $\hat{\G}$ satisfies (2). (Hence every 
homomorphism $\hat{\G}\to\outn$ has finite image.)
\end{remark}

\begin{remark} \label{r:ext} 
Let $\hat{G}$ be an arbitrary semisimple Lie group of real rank at least two, with centre $Z(\hat{G})$. Let $\hat{\G}<\hat{G}$ be an irreducible lattice, let $A=\hat{\G}\cap Z(\hat{G})$ and let $\G=\hat{\G}/A$. The abelianization of any subgroup of finite index in $\hat{\G}$ is finite (\cite{margulis}, page 333) and $\G$ is an irreducible lattice in the centreless semisimple Lie group $G=\hat{G}/Z(\hat{G})$. Thus the above remark allows us to remove from Corollary \ref{c:main} the hypothesis that the centre of $G$ is finite.
\end{remark}

In Theorem \ref{t:refinement}, condition (2) is used only to exclude the possibility
that a homomorphic image of $\G$ in $\outn$ might contain a fully irreducible
element. An alternative way of ruling out such images is to use bounded
cohomology, as in \cite{BF}.  We briefly review some notation.  A map $f:\G \to \mathbb{R}$ is a \emph{quasi-homomorphism} if the function $$(g,h) \mapsto |f(g) + f(h) - f(gh)|$$ is bounded on $\G \times \G$.  Let $V(\G)$ be the vector space of all quasi-homomorphisms from $\G$ to $\mathbb{R}$.  Two natural subspaces of $V(\G)$ are $B(\G)$, the vector space of bounded maps from $\G$ to $\mathbb{R}$, and $\hom(\G;\mathbb{R}),$ the vector space of genuine homomorphisms.   Define $\widetilde{\qh}(\G)=V(\G)/ (B(\G)+\hom(\G;\mathbb{R})).$

\begin{theorem} Let $\G$ be a group. Suppose that for every finite index subgroup $\L \subset \G$ the second
bounded cohomology of $\L$ is finite dimensional and  $\hom(\L,\mathbb{R})=0$.
Then every homomorphism $\phi:\G\to\outn$ has finite image.
\end{theorem}

\begin{proof}
Bestvina and Feighn \cite{BF} show that if $H < \out(F_n)$ contains a fully irreducible automorphism then either $H$ is virtually cyclic or $\widetilde{\qh}(H)$ is infinite dimensional.  If $\hom(\L;\mathbb{R})=0$ then a surjective map $\L\to H$ induces an injection $\widetilde{\qh}(H) \to \widetilde{\qh}(\L).$  The vector space $\widetilde{\qh}(\L)$ injects into the second bounded cohomology of $\L$ (see \cite{Micm}).  Therefore, for all finite index subgroups $\L \subset \G$ and integers $m$ the image of a homomorphism $\L\to\out(F_m)$ cannot contain a fully irreducible automorphism.  It follows from Corollary~\ref{c:BL} and the arguments in Proposition~\ref{p:IA_n} that $\phi(\G)$ is finite.
\end{proof}

In the light of \cite{BM}, this alternative to Theorem \ref{thm} also
implies Corollary \ref{c:main}.  

\bibliographystyle{plain}

\end{document}